\documentclass[12pt]{amsart}
\usepackage{amsmath,amssymb,amsbsy,amsfonts,latexsym,amsopn,amstext,cite,
                                               amsxtra,euscript,amscd,bm,mathabx}
\usepackage{url}
\usepackage[colorlinks,linkcolor=blue,anchorcolor=blue,citecolor=blue,backref=page]{hyperref}
\usepackage{color}
\usepackage{graphics,epsfig}
\usepackage{graphicx}
\usepackage{float} 
\usepackage[english]{babel}
\usepackage{mathtools}
\usepackage{todonotes}
\usepackage{url}
\usepackage[colorlinks,linkcolor=blue,anchorcolor=blue,citecolor=blue,backref=page]{hyperref}

\usepackage[norefs,nocites]{refcheck}

\hypersetup{breaklinks=true}

\usepackage[norefs,nocites]{refcheck}
\usepackage[english]{babel}
\begin{document}

\newtheorem{thm}{Theorem}
\newtheorem{lem}[thm]{Lemma}
\newtheorem{claim}[thm]{Claim}
\newtheorem{cor}[thm]{Corollary}
\newtheorem{prop}[thm]{Proposition} 
\newtheorem{definition}[thm]{Definition}
\newtheorem{question}[thm]{Open Question}
\newtheorem{conj}[thm]{Conjecture}
\newtheorem{prob}{Problem}

\newtheorem{lemma}[thm]{Lemma}

\newcommand{\hh}{{{\mathrm h}}}

\numberwithin{equation}{section}
\numberwithin{thm}{section}
\numberwithin{table}{section}

\def\vol {{\mathrm{vol\,}}}
\def\squareforqed{\hbox{\rlap{$\sqcap$}$\sqcup$}}
\def\qed{\ifmmode\squareforqed\else{\unskip\nobreak\hfil
\penalty50\hskip1em\null\nobreak\hfil\squareforqed
\parfillskip=0pt\finalhyphendemerits=0\endgraf}\fi}

\def \balpha{\bm{\alpha}}
\def \bbeta{\bm{\beta}}
\def \bgamma{\bm{\gamma}}
\def \blambda{\bm{\lambda}}
\def \bchi{\bm{\chi}}
\def \bphi{\bm{\varphi}}
\def \bpsi{\bm{\psi}}
\def \bomega{\bm{\omega}}
\def \btheta{\bm{\vartheta}}

\newcommand{\bfxi}{{\boldsymbol{\xi}}}
\newcommand{\bfrho}{{\boldsymbol{\rho}}}

\def\Kab{\cK_\psi(a,b)}
\def\Kuv{\cK_\psi(u,v)}
\def\SaUV{\cS_\psi(\balpha;\cU,\cV)}
\def\SaAV{\cS_\psi(\balpha;\cA,\cV)}

\def\SUV{\cS_\psi(\cU,\cV)}
\def\SAB{\cS_\psi(\cA,\cB)}

\def\Kmnp{\cK_p(m,n)}

\def\KKap{\cH_p(a)}
\def\KKaq{\cH_q(a)}
\def\KKmnp{\cH_p(m,n)}
\def\KKmnq{\cH_q(m,n)}

\def\Klmnp{\cK_p(\ell, m,n)}
\def\Klmnq{\cK_q(\ell, m,n)}

\def \SALMNq {\cS_q(\balpha;\cL,\cI,\cJ)}
\def \SALMNp {\cS_p(\balpha;\cL,\cI,\cJ)}

\def \SACXMQX {\fS(\balpha,\bzeta, \bxi; M,Q,X)}

\def\SAMJp{\cS_p(\balpha;\cM,\cJ)}
\def\SAMJq{\cS_q(\balpha;\cM,\cJ)}
\def\SAqMJq{\cS_q(\balpha_q;\cM,\cJ)}
\def\SAJq{\cS_q(\balpha;\cJ)}
\def\SAqJq{\cS_q(\balpha_q;\cJ)}
\def\SAIJp{\cS_p(\balpha;\cI,\cJ)}
\def\SAIJq{\cS_q(\balpha;\cI,\cJ)}

\def\RIJp{\cR_p(\cI,\cJ)}
\def\RIJq{\cR_q(\cI,\cJ)}

\def\TWXJp{\cT_p(\bomega;\cX,\cJ)}
\def\TWXJq{\cT_q(\bomega;\cX,\cJ)}
\def\TWpXJp{\cT_p(\bomega_p;\cX,\cJ)}
\def\TWqXJq{\cT_q(\bomega_q;\cX,\cJ)}
\def\TWJq{\cT_q(\bomega;\cJ)}
\def\TWqJq{\cT_q(\bomega_q;\cJ)}

 \def \xbar{\overline x}
  \def \ybar{\overline y}

\def\cA{{\mathcal A}}
\def\cB{{\mathcal B}}
\def\cC{{\mathcal C}}
\def\cD{{\mathcal D}}
\def\cE{{\mathcal E}}
\def\cF{{\mathcal F}}
\def\cG{{\mathcal G}}
\def\cH{{\mathcal H}}
\def\cI{{\mathcal I}}
\def\cJ{{\mathcal J}}
\def\cK{{\mathcal K}}
\def\cL{{\mathcal L}}
\def\cM{{\mathcal M}}
\def\cN{{\mathcal N}}
\def\cO{{\mathcal O}}
\def\cP{{\mathcal P}}
\def\cQ{{\mathcal Q}}
\def\cR{{\mathcal R}}
\def\cS{{\mathcal S}}
\def\cT{{\mathcal T}}
\def\cU{{\mathcal U}}
\def\cV{{\mathcal V}}
\def\cW{{\mathcal W}}
\def\cX{{\mathcal X}}
\def\cY{{\mathcal Y}}
\def\cZ{{\mathcal Z}}

\def\NmQR{N(m;Q,R)}
\def\VmQR{\cV(m;Q,R)}

\def\Xm{\cX_m}

\def \A {{\mathbb A}}
\def \B {{\mathbb A}}
\def \C {{\mathbb C}}
\def \F {{\mathbb F}}
\def \G {{\mathbb G}}
\def \L {{\mathbb L}}
\def \K {{\mathbb K}}
\def \Q {{\mathbb Q}}
\def \R {{\mathbb R}}
\def \Z {{\mathbb Z}}
\def \fS{\mathfrak S}

\def\e{{\mathbf{\,e}}}
\def\ep{{\mathbf{\,e}}_p}
\def\eq{{\mathbf{\,e}}_q}

\def\\{\cr}
\def\({\left(}
\def\){\right)}
\def\fl#1{\left\lfloor#1\right\rfloor}
\def\rf#1{\left\lceil#1\right\rceil}

\def\Tr{{\mathrm{Tr}}}
\def\Im{{\mathrm{Im}}}

\def \bFp {\overline \F_p}

\newcommand{\pfrac}[2]{{\left(\frac{#1}{#2}\right)}}

\def \Prob{{\mathrm {}}}
\def\e{\mathbf{e}}
\def\ep{{\mathbf{\,e}}_p}
\def\epp{{\mathbf{\,e}}_{p^2}}
\def\em{{\mathbf{\,e}}_m}

\def\Res{\mathrm{Res}}
\def\Orb{\mathrm{Orb}}

\def\vec#1{\mathbf{#1}}
\def\flp#1{{\left\langle#1\right\rangle}_p}

\def\mand{\qquad\mbox{and}\qquad}
%
%
%
%


\title[Double Sums of Kloosterman sums in Finite Fields]
{Double Sums of Kloosterman Sums in Finite Fields}

  \author[S.  Macourt] {Simon Macourt}
\address{Department of Pure Mathematics, University of New South Wales,
Sydney, NSW 2052, Australia}
\email{s.macourt@unsw.edu.au}

\author[I. E. Shparlinski]{Igor E. Shparlinski} 
\address{School of Mathematics and Statistics, University of New South Wales, 
Sydney, NSW 2052, Australia}
\email{igor.shparlinski@unsw.edu.au}

\date{\today}

\begin{abstract}  We bound double sums of  Kloosterman sums over a  finite 
field $\F_{q}$, with one or both  parameters ranging over an affine space over 
its prime  subfield $\F_p \subseteq \F_{q} $. These are finite fields analogues
of a series of recent results by various authors in finite fields and residue rings.    
Our results are based on recent advances in additive combinatorics in arbitrary finite field. 
\end{abstract}



\maketitle

\section{Introduction}

 \subsection{Background}
 
 The motivation behind this work comes from recent advances in estimating 
 various bilinear sums of Kloosterman sums which have found a wealth of 
 applications to various arithmetic problems, see~\cite{BFKMM1,FKM,KMS1,KMS2,LSZ,Shp,ShpZha} and references therein. 

Here we extend some of these results to the settings of finite fields. Our approach 
is modelled from that of~\cite{Shp,ShpZha} however at one significant point it deviates and we employ some very recent results of~\cite{Moh} from additive combinatorics in 
arbitrary finite fields. 

For a prime power  $q$,  let $\F_q$ denote the   finite field $\K = \F_q$
of $q$ elements. 

We fix a nontrivial additive character $\psi$ of $\F_q$ and 
for integers $u,v \in \F_q$ we define
the {\it Kloosterman sum\/}
$$
\Kuv = \sum_{x \in \F_q^*} \psi\(ux + vx^{-1} \).
$$
We consider  sums of Kloosterman sums
$$
\SUV = \sum_{u\in \cU} \sum_{v \in \cV} \Kuv 
$$
over some subsets  $\cU, \cV \subseteq \F_q$ and also 
of more general sums, 
$$
\SaUV = \sum_{u\in \cU} \sum_{v \in \cV} \alpha_v \Kuv,
$$
with a   sequence of complex weights $\balpha = \{\alpha_v\}_{v\in \cV}$.

By the Weil bound we have
$$
|\Kuv|\le  2 q^{1/2},
$$
see~\cite[Corollary~11.12]{IwKow}. Hence we immediately obtain 
 \begin{equation}
\label{eq:trivial AUV}
\left| \SaUV \right| \le 2  U q^{1/2} \sum_{v\in \cV}  |\alpha_v|,
\end{equation}
where $U = \# \cU$ is the cardinality of $\cU$.

We are interested in studying cancellations amongst Kloosterman sums and
thus improvements of the trivial bound~\eqref{eq:trivial AUV}. We note that 
the sums $\SUV$ and $\SaUV$ are finite field analogues of similar sums 
studied in~\cite{BFKMM1,FKM,KMS1,KMS2,LSZ,Shp,ShpZha} in the settings 
of prime fields and residue rings. Hence, adopting the model of $\F_q$ as 
$$
\F_q \cong \F_p[X]/f(x)
$$
for a prime field $\F_p$ and an irreducible over $\F_p$ polynomial $f$ (of degree 
$\deg f = [\F_q:\F_p]$)  we expect that our bounds can be used
for similar arithmetic applications in function fields. Since in the above works, 
the case of averaging over intervals plays a vital role, here we
consider the case when one or both of the sets $\cU$ and $\cV$ is an affine 
subspace of $\F_q$, considered as a vector space over its prime subfield 
$\F_p \subseteq \F_{q}$. More precisely, we consider the sums
\begin{itemize}
\item  $\SAB$  with two affine spaces $\cA$ and $\cB$, 
\item $\SaAV$  with an affine space  $\cA$ and an arbitrary set  $\cV$. 
\end{itemize}
Our approach is similar to that of~\cite{Shp,ShpZha} however some ingredients 
used in~\cite{Shp,ShpZha} are either unknown or do not exist in function field settings.
Hence we use a different approach based on additive combinatorics and in particular we
rely on recent results of Mohammadi~\cite{Moh}.

 \subsection{Our results}

We recall that the
notations $U = O(V)$, $U \ll V$ and  $V \gg U$ are all equivalent to the
statement that the inequality $|U| \le c V$ holds with some
constant $c> 0$, which is absolute throughout this paper.

We start with the sums $\SAB$. 

\begin{thm}
\label{thm:SAB}  Assume that $q = p^{2k+1}$ is an odd power of a prime $p$. 
Let $\cA$ and $\cB$ be affine subspaces of $\F_q$ of cardinalities $A=\#\cA$ 
and $B = \# \cB$, respectively, with. $A\le B$.  Then  
$$ 
\SAB \ll AB \max\{ q^{52/153}, (q/A)^{831/832},  (q/A)^{761/760} q^{-1/760}\}.
$$
\end{thm}

Clearly Theorem~\eqref{thm:SAB} is only non-trivial  for  $A \ge q^{415/831}$ as otherwise the
bound $|\SAB| \le 2AB q^{1/2}$, implied by~\eqref{eq:trivial AUV}, is stronger.

Given a   sequence of complex weights $\balpha = \{\alpha_v\}_{v\in \cV}$ supported on 
a set $\cV$ and $\rho > 0$, as usual, we define
$$
\| \balpha\|_\rho = \(\sum_{v \in \cV} |\alpha_v|^\rho\)^{1/\rho}. 
$$
 
 \begin{thm}
\label{thm:SaAV}  Assume that $q = p^{2k+1}$ is an odd power of a prime $p$. 
Let $\cA$  be an affine subspace of $\F_q$ of cardinality $A=\#\cA$ 
and let $\cV \subseteq \F_q$ be an arbitrary set   of cardinality $V=\#\cV$.  Then, 
for any  sequence of complex weights $\balpha = \{\alpha_v\}_{v\in \cV}$ we have
\begin{align*}
\SaAV \ll  Aq^{1/4}  & \sqrt{|\balpha|_1 |\balpha|_2} \\&
\max\left\{q^{13/51}, \left(\frac{q}{A}\right)^{935/1248} , q^{-1/1140} \left( \frac{q}{A}\right)^{214/285}\right\}.
\end{align*}
\end{thm}
 
If we suppose that $|\alpha_v| \ll 1$ for all $v \in \cV$, then clearly Theorem~\ref{thm:SaAV} is only non-trivial for $A^{935/623}V^{312/623} > q$. If we suppose $V=A$ then we have Theorem~\ref{thm:SaAV} is non-trivial provided $A>q^{623/1247}$.

 \section{Background from additive combinatorics}

For a set $\cS  \subseteq \F_q$, we use $E(\cS)$ to denote its {\it additive 
energy\/}, that is, the number of solutions to the equation 
$$
s_1+s_2 =s_3 +s_4 , \qquad s_1, s_2, s_3,s_4\in \cS.
$$
Also, as usual, we denote
$$
\cS^{-1} = \{s^{-1}~:~s \in \cS\} \mand 
2 \cS =  \{s+t~:~s+t \in \cS\} .
$$

Then by~\cite[Corollary~5]{Moh} we have

\begin{lem}
\label{lem:E_1/S}  Let $\cS  \subseteq \F_q$  with 
$$
\# \cS = S \mand \#(2\cS) = T
$$
and  such that 
$T \ge (\# G)^{52/51}$ for any proper subfield $\G$ of $\F_q$.  Then 
$$ 
E\(\cS^{-1}\)   \ll   \(T^{173/104} + q^{-1/285} T^{476/285}\) S^{4/3}.
$$
\end{lem}

It is easy to see that the Cauchy inequality implies the well-known inequality
$$
S^4 \le \#\(2\cS^{-1}\)  E\(\cS^{-1}\) . 
$$
Hence from Lemma~\ref{lem:E_1/S} we derive the following (see also~\cite[Corollary~5]{Moh}).

\begin{cor}
\label{cor:2S_2/S}  Let $\cS  \subseteq \F_q$  with 
$$
\# \cS = S, \qquad \#(2\cS) = T, \qquad  \#(2\cS^{-1}) = U
$$
and  such that 
$T \ge (\# G)^{52/51}$ for any proper subfield $\G$ of $\F_q$.  Then 
$$ 
\max\{T,U\} \gg  \min\{S^{832/831}, q^{1/761} S^{760/761}\}.
$$
\end{cor}

\section{Proof of Theorem~\ref{thm:SAB}}

Changing the order of summation,  we obtain
 \begin{equation}
\label{eq:SAB Transf}
 \SAB = \sum_{x\in \F^*_q}
  \sum_{a\in \cA}\psi(ax) \sum_{b \in \cB}\psi(bx^{-1} ).
\end{equation}

Clearly if $\cA$ is a translate of a linear space $\cL$ then 
 \begin{equation}
\label{eq:Orth A}
 \sum_{a\in \cA}\psi(ax) =
 \begin{cases} A, & \text{if}\  x\in \cL^\perp, \\
0, & \text{otherwise},
 \end{cases} 
\end{equation}
where $\cL^\perp$ denotes the the orthogonal complement to $\cL$.

Similarly  if $\cB$ is a translate of a linear space $\cM$ then 
 \begin{equation}
\label{eq:Orth B}
 \sum_{b\in \cB}\psi(bx^{-1}) =
 \begin{cases} B, & \text{if}\  x^{-1}\in \cM^\perp, \\
0, & \text{otherwise}.
 \end{cases} 
\end{equation}

Hence, substituting~\eqref{eq:Orth A} and~\eqref{eq:Orth B} in~\eqref{eq:SAB Transf}, 
we obtain
 \begin{equation}
\label{eq:SABS}
 \SAB \ll AB S
\end{equation}
 where $S = \#\cS$ and the set $\cS$ is defined as follows
 $$
 \cS = \{x \in \F_q^*~:~x\in \cL^\perp\ \text{and}\  x^{-1}\in \cM^\perp\}.
 $$

 If $S \le q^{52/153}$ the result follows immediately.
 Otherwise we see that
 $$
 \#(2\cS) \ge \#\cS \ge   (\# \G)^{52/51}
$$ 
for any proper subfield $\G$ of $\F_q$ (since $q$ is not a perfect square
we have $\# \G\le q^{1/3}$). Hence Corollary~\ref{cor:2S_2/S} applies to $\cS$.

Since $\cL^\perp$ and $\cM^\perp$ are  linear spaces, we  obviously have
$$
2 \cS \subseteq \cL^\perp \mand 2 \cS^{-1} \subseteq \cM^\perp.
$$
Consequently, 
$$
\# (2 \cS) \le  \# \cL^\perp = q/A \mand 
\# \(2 \cS^{-1}\) \le  \# \cM^\perp = q/B. 
$$
Invoking  Corollary~\ref{cor:2S_2/S}, we obtain
$$
 \min\{S^{832/831}, q^{1/761} S^{760/761}\} \ll \max\{q/A, q/B\} \le q/A
$$
as $A \le B$. Thus
$$
S \ll \max\{(q/A)^{831/832},  (q/A)^{761/760} q^{-1/760}\}
$$
which after substitution in~\eqref{eq:SABS} implies the result. 

\section{Proof of Theorem~\ref{thm:SaAV}}

By changing the order of summation we have
\begin{align*}
|\SaAV| = \left|\sum_{x \in \F^*_q}\sum_{v \in \cV} \alpha_v \psi(vx^{-1}) \sum_{a \in \cA} \psi(ax)\right|.
\end{align*}
As previously, if $\cA$ is a translate of a linear space $\cL$ then
\begin{equation}
\label{eq:Orth A2}
 \sum_{a\in \cA}\psi(ax) =
 \begin{cases} A, & \text{if}\  x\in \cL^\perp, \\
0, & \text{otherwise,}
 \end{cases} 
\end{equation}
where $\cL^\perp$ denotes the orthogonal complement $\cL$. It follows that
\begin{align*}
|\SaAV| &=A \left|\sum_{x \in \cL^\perp}\sum_{v \in \cV} \alpha_v \psi(vx^{-1})\right| \\
&\le A \left|\sum_{v \in \cV} \alpha_v\right| \left| \sum_{x \in \cL^\perp} \psi(vx^{-1})\right| \\
&\le A \sum_{v \in \cV} \left|\alpha_v\right|^{1/2}\left|\alpha_v^2\right|^\frac{1}{4}  \left| \sum_{x \in \cL^\perp} \psi(vx^{-1})\right|.
\end{align*}
Applying the Cauchy--Schwartz inequality twice, we obtain
\begin{align*}
|\SaAV|^4 &\le A^4 \left(\sum_{v \in \cV} \left|\alpha_v\right| \right)^{2} \sum_{v \in \cV}|\alpha_v|^2\sum_{v \in \cV} \left|\sum_{x \in \cL^\perp} \psi(vx^{-1})\right|^4  \\
&\le A^4|\balpha|_1^2 |\balpha|_2^2\\
& \qquad \quad \sum_{v \in \F^*_q} \sum_{w,x,y,z \in \cL^\perp} \psi(v(w^{-1}+x^{-1}-y^{-1}-z^{-1})) .
\end{align*}

If $2\cL^\perp=\cL^\perp < (\#G)^{52/51}$, we use the trivial bound on additive energy, that is $E((\cL^\perp)^{-1}) \le (\# \cL^\perp)^3$, and the result follows immediately. Otherwise, we apply Lemma~\ref{lem:E_1/S}, observing $\#L^\perp=q/A$, to obtain 
\begin{align*}
\SaAV &\ll Aq^{1/4} \sqrt{|\balpha|_1 |\balpha|_2}\\
& \qquad \quad   \left(\left(\left(\frac{q}{A}\right)^{173/104}+q^{-1/285}\left(\frac{q}{A}\right)^{476/285}\right)\left(\frac{q}{A}\right)^{4/3}\right)^{1/4} \\
&\ll Aq^{1/4}  \sqrt{|\balpha|_1 |\balpha|_2} \\
& \qquad \quad  \max\left\{\left(\frac{q}{A}\right)^{935/1248} , q^{-1/1140} \left( \frac{q}{A}\right)^{214/285}\right\}.
\end{align*}
This completes the proof.

\section{Comments}

Some of the motivation to this paper comes from an intention to obtain function field analogues of the 
asymptotic formulas, with a power saving, from~\cite{BFKMM1, BFKMM2, Shp, Young} 
for 4th moments of $L$-functions. However, despite recent progress in this direction due 
to Florea~\cite{Flo},  some ingredients, used in the groundbreaking work of 
Young~\cite{Young},  remain missing in the function field case. 

Finally, we need to impose the condition on   $q = p^{2k+1}$ to avoid th existence of large subfields.
It is certainly interesting to drop this restriction and extend our results to even degree extensions. 

\section*{Acknowledgements}

The authors are grateful to   Goran Djankovi{\'c} for very illuminating discussions.

This work was  supported in part  by ARC Grant~DP170100786.

\end{document}